\documentclass{amsart}
\title{Generating functions for the number of curves on Abelian surfaces}
\author{Jim Bryan\\
Naichung Conan Leung}
\date{\today}
\address{Department of Mathematics\\
University of California\\
Berkeley, CA 94720}
\address{School of Mathematics\\
University of Minnesota\\
Minneapolis, MN 55455}

\usepackage{amscd,eepic,epic}

\newtheorem{thm}{Theorem}[section]
\newtheorem{theorem}[thm]{Theorem}
\newtheorem{cor}[thm]{Corollary}

\newtheorem{lemma}[thm]{Lemma}
\newtheorem{prop}[thm]{Proposition}

\newtheorem{definition}[thm]{Definition}

\newtheorem{rem1}[thm]{Remark}

\newtheorem{exer1}[thm]{Exercise}

\newcommand{\cnums} {{\mathbf C}}          
\newcommand{\rnums} {{\mathbf R}}		
\newcommand{\znums} {{\mathbf Z}}		

\newcommand{\dbar}{\bar\partial}
\newcommand{\del}{\partial}

\newcommand{\im}{\operatorname{Im}}

\newcommand{\ie}{{\em i.e. }}

\newcommand{\point}{\text{pt.}}
\newcommand{\til}[1]{\Tilde{#1}}

\renewcommand{\P}{\mathbf{P}}
\newcommand{\M}{\mathcal{M}}

\newcommand{\pic}{\operatorname{Pic}}
\newcommand{\pf}{{\sc Proof: }}
\newcommand{\virdim}{\operatorname{virdim}}

\begin{document}
\begin{abstract}
Let $X$ be an Abelian surface and $C$ a holomorphic curve in $X$ representing a
primitive homology class. The space of  genus $g$  curves in the class of
$C$ is $g$ dimensional. We count the number of such curves that pass
through $g$ generic points  and we also count the number of curves in the
fixed linear system $|C|$ passing through $g-2$ generic points. These two
numbers, (defined appropriately) only depend on $n$ and $g$ where
$n=\frac{C\cdot C}{2} +1-g$ and not on the particular $X$ or $C$ ($n$ is
the number of nodes when a curve is nodal and reduced). 

G\"ottsche conjectured that certain quasi-modular forms are the  generating
functions for the number of curves in a fixed linear system
\cite{Gott-conj}. Our theorem proves his 
formulas and shows that (a different) modular form also arises in the problem of
counting curves without fixing a linear system.  We use techniques that 
were developed in \cite{Br-Le} for similar questions on $K3$ surfaces. The
techniques include Gromov-Witten invariants for families and a degeneration
to an elliptic fibration. One new feature of the Abelian surface case is
the presence of 
non-trivial $\pic^{0} (X)$. We show that for any surface $S$ the cycle in the
moduli space of stable maps defined by requiring that the image of the map
lies in a fixed linear system is homologous to the cycle defined by
requiring the image of the map meets $b_{1}$ generic loops in $S$
representing the generators of $H_{1}(S;\znums )/\operatorname{Tor}$.  
\end{abstract}
\thanks{The first author is supported by a grant from the Ford Foundation
and the second author is supported by NSF grant DMS-9626689.} 
\maketitle
\markboth{Curves on Abelian surfaces}
{Curves on Abelian surfaces}
\renewcommand{\sectionmark}[1]{}

\tableofcontents
\pagebreak

\section{Introduction}\label{sec: intro}
Let $X$ be an Abelian surface and let $C$ be a holomorphic curve in $X$
representing a primitive homology class. For $n$ and $g$ satisfying $C\cdot
C=2g-2+2n$,  there is a $g$ dimensional space of curves of genus $g$ in the
class of $C$. To define an enumerative problem, one must impose $g$
constraints on the curves. There are two natural ways to do this. One way
is to count the number of  curves passing through $g$ generic points which
we denote $N_{g,n}(X,C)$. The second way is to count the number of curves
in the fixed linear system $|C|$ 
passing through $g-2$ generic points which we denote $N_{g,n}^{FLS}(X,C)$.  We
define (modified) Gromov-Witten invariants that compute the numbers
$N_{g,n}(X,C)$ and $N_{g,n}^{FLS}(X,C)$ 
and we prove that they do not depend on $X$ or $C$ but are universal
numbers henceforth denoted $N_{g,n}$ and $N^{FLS}_{g,n}$. Our main theorem
computes these numbers  as the Fourier coefficients  of quasi-modular
forms. Note that $X$ does not contain any genus 0 curves
so implicitly $g>0$ throughout. 

\begin{thm}[Main theorem]\label{thm: main thm}
The universal numbers $N^{FLS}_{g,n}$ and $N_{g,n}$ are given by the
following generating functions:
\begin{eqnarray*}
\sum_{n=0}^{\infty }N_{g,n}q^{n+g-1}&=&g(DG_{2})^{g-1},\\
\sum_{n=0}^{\infty }N_{g,n}^{FLS}q^{n+g-1}&=&(DG_{2})^{g-2}D^{2}G_{2}\\
&=&(g-1)^{-1}D((DG_{2})^{g-1})
\end{eqnarray*}
where $D$ is the operator $q\frac{d}{dq}$ and $G_{2}$ is the Eisenstein
series, \ie 
$$
G_{2}(q)=-\frac{1}{24}+\sum_{k=1}^{\infty }(\sum_{d|k}d)q^{k}. 
$$
\end{thm}

Note that the right hand sides are
quasi-modular forms (quasi-modular forms are an algebra generated by
modular forms, $G_{2}$, and $D$, c.f. \cite{Gott-conj}) and the left hand
sides are reminiscent of theta series since the power of $q$ is
$\frac{C\cdot C}{2}=n+g-1$.

The formula for $N_{g,n }^{FLS}$ in Theorem \ref{thm: main thm} was
first conjectured by Lothar 
G\"ottsche \cite{Gott-conj} who proved the genus 2 case (see also
\cite{DeB}). More generally, G\"ottsche gave conjectural generating
functions for 
the number of curves with $n$ nodes in any $n$-dimensional sublinear system of
$|C|$ for any sufficiently ample divisor $C$ on any surface $S$. His
formulas involve only 
$c_{2}(S)$,  $C\cdot C$, $C\cdot K$, $K\cdot K$ and universal functions
($K$ is the canonical class). For surfaces with numerically trivial
canonical class, his formulas reduce to quasi-modular forms.

For concreteness, we expand the equations of the theorem to write:
\begin{eqnarray*}
\sum_{n=0}^{\infty }N_{g,n}q^{n+g-1}&=&g\left(\sum_{k=1}^{\infty
}k(\sum_{d|k}d)q^{k} \right)^{g-1},\\ 
\sum_{n=0}^{\infty }N_{g,n }^{FLS}q^{n+g-1}&=&\left(\sum_{k=1}^{\infty
}k(\sum_{d|k}d)q^{k} 
\right)^{g-2}\left(\sum_{k=1}^{\infty }k^{2}(\sum_{d|k}d)q^{k} \right).\\ 
\end{eqnarray*}

So for example, for small values of $n $ and $g$,  $N_{g,n}^{FLS}$
and $N_{g,n}$ are given in the following tables: 
\vskip .25in
\begin{tabular}{|c||c|c|c|c|c|c|c|c|}\hline 
$N_{g,n}^{FLS}$&    n=	0&     n=1&	n=2&	n=3&	n=4&	n=5&
n=6&	n=7\\ \hline \hline  
g=2   &	1&	12&	36&	112&	150&	432&	392&	960\\ \hline 
g=3   &	1&	18&	120&	500&	1620&	4116&	9920&	19440\\ \hline 
g=4   &	1&	24&	240&	1464&	6594&	23808&	73008&	198480\\ \hline
g=5   &	1&	30&	396&	3220&	18960&	88452&	344960&	1169520 \\ \hline 	
\end{tabular}
\vskip .25in
\begin{tabular}{|c||c|c|c|c|c|c|c|c|}\hline 
$N_{g,n}$&    n=	0&     n=1&	n=2&	n=3&	n=4&	n=5&
n=6&	n=7\\ \hline \hline  
g=2  & 2&    12&    24&        56&      60&      144&      112&
240\\ \hline 
g=3&   3&    36&    180&    600&    1620&     3528&     7440&       12960\\
\hline 
g=4&   4&    72&    576&    2928&   11304&   35712&   97344&     238176\\
\hline  
g=5&   5&    120&   1320&  9200&  47400&  196560&   689920&   2126400\\
\hline  
\end{tabular}

\vskip .25in

Gromov-Witten invariants have been remarkably effective in answering many
questions in enumerative geometry for certain varieties such as $\P ^{n}$;
however, the ordinary Gromov-Witten invariants are not very effective for
counting 
curves on most surfaces. One basic reason is that the moduli space of
stable maps often fails to be a good model for a linear system (and the
corresponding Severi varieties) for dimensional
reasons. For an bundle $L$ such that $L-K$ is ample, the dimension of the
Severi variety 
$V_{g}(L)$ (the closure of the set of geometric  genus $g$ curves 
in the  complete linear  system $|L|$) is 
$$
\dim _{\cnums }V_{g}(L)=-K\cdot L +g-1+p_{g}-q
$$
where $K$ is the canonical class, $p_{g}=\dim H^{0}(X,K)$, and $q=\dim
H^{1}(X,\mathcal{O})$. On the other hand, the virtual dimension of the
moduli space $\M _{g,L}(X)$ of stable maps of genus $g$ in the class
dual to $c_{1}(L)$ is 
$$
\virdim _{\cnums }\M _{g,L}(X)=-K\cdot L+g-1.
$$

The discrepancy $p_{g}-q$ arises from two sources. Since the image of maps
in $\M _{g,L}(X)$ are divisors not only in $|L|$ but also potentially in
every linear 
system in $\pic ^{c_{1}(L)}(X)$,  one would expect $\dim \M _{g,L}(X)$ to
exceed $\dim V_{g}(L)$ by $q=\dim \pic ^{c_{1}(L)}(X)$ (we use $\pic
^{c}(X)$ to denote the component of $\pic (X)$ with Chern class
$c$). We show that this 
discrepancy can be accounted for within the framework of the usual
Gromov-Witten invariants (see Theorem \ref{thm: pullback of class on pic is
ft*ev*}). 

However, even if we consider $\M _{g,L}$ as a model for the parametrized
Severi varieties
$$
V_{g}(c_{1}(L))\equiv \bigcup _{L'\in \pic^{c_{1}(L)}(X)}V_{g}(L'), 
$$
there is still a $p_{g}$ dimensional discrepancy (see also \cite{Do}).

The reason is the following.  The virtual dimension
of $\M _{g,L}(X) $ is the dimension of the space of curves that persist as
{\em pseudo-holomorphic} curves when we perturb the K\"ahler structure to a
generic almost K\"ahler structure. The difference of $p_{g}$ in the
dimensions of $\M _{g,L}$ and $V_{g}(c_{1}(L))$ means that only a
codimension $p_{g}$ subspace  of $V_{g}(c_{1}(L))$ persists as
pseudo-holomorphic curves when we perturb the 
K\"ahler structure. One way to rectify this situation is to find a compact
$p_{g}$-dimensional\footnote{By this we mean a real $2p_{g}$
dimensional family. Note that the parameter space for the family need not
have an 
almost complex structure.}  {\em family} of almost K\"ahler structures that
has the 
property that the only almost K\"ahler structure in the family that
supports pseudo-holomorphic curves in the class $c_{1}(L)$ is the original
K\"ahler structure. If $T$ is such a family, then the moduli space $\M
_{g,L}(X,T)$ of stable maps for the family $T$ is a better model for the
space $V_{g}(c_{1}(L))$ in the sense that its dimension is stable under
generic perturbations of the family $T\mapsto T'$.

It is straight-forward to extend the notion of the ordinary Gromov-Witten
invariants to invariants
for families of almost K\"ahler structures. The invariants will only depend
on the deformation class of the underlying family of symplectic structures.
Given the existence of a $p_{g}$-dimensional family as described above,
these invariants can be used to answer enumerative geometry questions for
the corresponding surface and linear system. 

In general, it is not clear when such a family will exist; however,
if $X$ has a hyperk\"ahler metric $g$ (\ie $X$ is an Abelian or $K3$ surface),
then there is a natural candidate for $T$, namely the hyperk\"ahler family
of K\"ahler structures. We call this family the {\em twistor family}
associated to the metric  $g$ and we denote it $T_{g}$. It is parameterized
by a 2-sphere and so $\dim _{\rnums }T_{g}=2=2p_{g}$ as it should.
Furthermore, the property that all the curves in $\M _{g,L}(X,T_{g})$ are
holomorphic for the original complex structure can be proved with Hodge
theory (of course this need no longer be the case for a  perturbation of 
$T_{g}$ to a generic family of almost K\"ahler structures).

We will define the numbers $N_{g,n}^{FLS}$ and $N_{g,n}$ as certain
Gromov-Witten invariants
for the twistor family $T_{g}$ associated to a hyperk\"ahler metric on an
Abelian surface $X$ (see our previous paper \cite{Br-Le} for the $K3$
case). We show that the invariants only depend on $g$ and $n$
(not $X$ or $C$) and they count each irreducible geometric genus $g$ curve
with positive integral multiplicity. Furthermore, the multiplicity is 1 if
additionally the curve is nodal.   There are additional Gromov-Witten
invariants for the twistor family that we also compute. These invariants
are easier to describe with the notation of the Gromov-Witten invariants so
we will postpone the statement of the result until section \ref{sec:
twistor family} (Defintion \ref{def: definition of Ngn and NngFLS} and
Theorem \ref{thm: restatement of main thm including Nij}).   The
enumerative problem these additional invariants correspond to is counting
curves that pass through $g-1$ points and lie in a certain one dimensional
family of linear systems.

Section \ref{sec: GW for families}  reviews Gromov-Witten invariants for
families and formulates 
our result that equates the subset of $\M _{g,L}(X)$ consisting of maps
whose image have fixed divisor class with a cycle more familiar in ordinary
Gromov-Witten theory (see Theorem \ref{thm: pullback of class on pic is
ft*ev*}). In section \ref{sec: twistor family}
we discuss properties
of the twistor family, define $N_{g,n}^{FLS}$ and $N_{g,n}$, and prove they
have the enumerative 
properties discussed above. In section \ref{sec: computing Ng(n)}  we
compute the invariants to 
complete the proof of our main theorem and its generalization. We conclude
with an appendix 
containing the proof of Theorem \ref{thm: pullback of class on pic is
ft*ev*} that was postponed in the main exposition.

The authors would like to thank O. DeBarre, A Givental, L. G\"ottsche, K.
Kedlaya, and C. Taubes for helpful discussions and correspondence.

\section{Gromov-Witten invariants for families}\label{sec: GW for families}
We review Gromov-Witten invariants for families and we refer the reader to
\cite{Br-Le} for details. 

Let $X$ be any compact symplectic manifold with an almost K\"ahler
structure. 
Recall that an $n$-marked, genus $g$ {\em stable map} of degree $C\in
H_{2}(X,\znums )$ is a (pseudo)-holomorphic map $f:\Sigma \to X$ from an
$n$-marked nodal curve $(\Sigma ,x_{1},\ldots,x_{n})$ of geometric genus
$g$ to $X$ with $f_{*}([\Sigma ])=C$ that has no infinitesimal
automorphisms. Two 
stable maps $f:\Sigma \to  X$ and $f':\Sigma' \to X $ are equivalent if
there is a biholomorphism $h:\Sigma \to \Sigma '$ such that $f=f'\circ h$.
We write $\M _{g,n,C}(X,\omega )$ for the moduli space of equivalence
classes of genus $g$, $n$-marked, stable maps of degree $C$ to $X$. We will
often drop the $\omega $ or $X$ from the notation if they are understood
and we sometimes will drop the $n$ from the notation when it is $0$. If $B$
is a family of almost K\"ahler structures, we denote parameterized version
of the moduli space:
$$
\M _{g,n,C}(X,B)=\bigcup_{t\in B}\M _{g,n,C}(X,\omega _{t}).
$$

If $B$ is a compact, connected, oriented manifold then $\M _{g,n,C}(X,B)$
has a fiduciary cycle $[\M _{g,n,C}(X,B)]^{vir}$ called the virtual
fundamental cycle (see \cite{Br-Le} and the fundamental papers of Li and
Tian \cite{Li-Tian2}\cite{Li-Tian}\cite{Li-Tian3}). The dimension of the
cycle is
$$
\dim_{\rnums }[\M _{g,n,C}(X,B)]^{vir}=2c_{1}(X)(C)+(6-\dim _{\rnums
}X)(g-1)+ 2n+\dim _{\rnums }B.
$$

The invariants are defined by evaluating cohomology classes of $\M _{g,n,C}$
on the virtual fundamental cycle. The cohomology classes are defined via
incidence relations of the maps with cycles in $X$. The framework is as
follows. There are maps
$$
\begin{CD}
\M _{g,1,C}@>{ev}>>X\\
@VV{ft}V\\
\M _{g,C}
\end{CD}
$$
called the {\em evaluation} and {\em forgetful} maps defined by
$ev(\{f:(\Sigma ,x_{1})\to X \})=f(x_{1})$ and $ft(\{f:(\Sigma ,x_{1})\to X
\})=\{f:\Sigma \to X \}$.\footnote{There is some subtlety to making this
definition rigorous since forgetting the point may make a stable map
unstable, but it can be done.} The diagram should be regarded as the
universal map over $\M _{g,C}$.  

Given geometric cycles $\alpha _{1},\ldots,\alpha _{l}$ in $X$ representing
classes $[\alpha _{1}],\ldots,[\alpha _{l}]\in H_{*}(X,\znums )$ with
Poincar\'e duals $[\alpha _{1}]^{\vee},\ldots,[\alpha _{l}]^{\vee}$, we can
define the Gromov-Witten invariant 
\begin{equation*}
\Phi _{g,C}^{(X,B)}(\alpha _{1},\ldots,\alpha _{l})=\int_{[\M
_{g,C}(X,B)]^{vir}} ft_{*}ev^{*}([\alpha _{1}]^{\vee })\cup \cdots \cup
ft_{*}ev^{*}([\alpha _{l}]^{\vee }) .
\end{equation*}
$\Phi _{g,C}^{(X,B)}(\alpha _{1},\ldots,\alpha _{l})$ counts the number of
genus $g$, degree $C$ maps which are pseudo-holomorphic with respect to
some almost K\"ahler structure in $B$ and such that the image of the map
intersects each of the cycles $\alpha _{1},\ldots,\alpha
_{l}$.\footnote{The integral is defined to be 0 if the integrand is not a
class of the correct degree.} The Gromov-Witten invariants are multi-linear
in the $\alpha $'s and they are symmetric for $\alpha $'s of even degree
and skew symmetric for $\alpha $'s of odd degree. If $p_{1},\ldots,p_{k }$
are points in  a path-connected $X$, we will use the shorthand
$$
\Phi^{(X,B)}_{g,C}(\point ^{k},\alpha _{k+1},\ldots,\alpha _{l}):=\Phi
^{(X,B)}_{g,C}(p_{1},\ldots,p_{k},\alpha _{k+1},\ldots,\alpha_{l} ). 
$$

Now suppose that $X$ is a K\"ahler surface.
In order to count curves in a fixed linear system $|L|$ with
$c_{1}(L)=[C]^{\vee }$ one would like to restrict the above integral 
 to the cycle defined by $\Psi_{\Sigma _{0}} ^{-1}(0)$ where $\Psi_{\Sigma
_{0}} $ is the map
$$
\Psi_{\Sigma _{0}} :\M _{g,C}(X,\omega )\to \pic ^{0}(X)
$$
given by $f\mapsto \mathcal{O}(\im (f)-\Sigma _{0})$ where $\Sigma _{0}\in
|L|$ is a fixed divisor. Dually, one can add the pullback by
$\Psi_{\Sigma _{0}} $ of the volume form on $\pic ^{0}(X)$ to the integrand
defining the invariant: 
$$
\int_{[\M _{g,C}(X)]^{vir}}\Psi_{\Sigma _{0}} ^{*}([\point ]^{\vee })\cup
ft_{*}ev^{*}([\alpha 
_{1}]^{\vee }) \cup \cdots \cup ft_{*}ev^{*}([\alpha _{l}]^{\vee }).
$$

We show that $\Psi_{\Sigma _{0}} ^{*}([\point ]^{\vee })$ can be expressed
in the usual Gromov-Witten framework.
\begin{thm}\label{thm: pullback of class on pic is ft*ev*}
Let $X$ be a K\"ahler surface and let $[\gamma]\in H_{1}(X,\znums )  $ and
let $\til{\gamma }$ be the 
corresponding class in $H^{1}(\pic ^{0}(X),\znums )$ induced by the
identification $\pic ^{0}(X)\cong H^{1}(X,\rnums )/H^{1}(X,\znums )$. Then
$$
\Psi _{\Sigma _{0}}^{*}(\til{\gamma})=ft_{*}ev^{*}([{\gamma}]^{\vee }). 
$$ 
\end{thm}
\begin{cor}
Let $[\gamma _{1}],\ldots,[\gamma _{b_{1}}]$ be an oriented integral basis
for $H_{1}(X;\znums )$. Then 
$\Psi_{\Sigma _{0}} ^{*}([\point ]^{\vee })=ft_{*}ev^{*}([\gamma
_{1}]^{\vee })\cup \cdots \cup 
ft_{*}ev^{*}([\gamma _{b_{1}}]^{\vee }).$ 
\end{cor}
\pf We defer the proof of Theorem \ref{thm: pullback of class on pic is
ft*ev*} to the appendix. The corollary follows immediately.

The upshot is that we count curves in a fixed linear system using the usual
constraints from Gromov-Witten theory. Namely, the invariant:
$$
\Phi _{g,C}(\gamma _{1},\ldots,\gamma _{b_{1}},\alpha _{1},\ldots,\alpha _{l})
$$
counts the number of genus $g$ maps whose image lie in a fixed linear
system $|L|$ with $c_{1}(L)=[C]^{\vee }$ and  hit the cycles $\alpha
_{1},\ldots,\alpha _{l}$. 

\section{The twistor family and the definition of 
$N_{g,n}$ and $N^{FLS}_{g,n}$}\label{sec: twistor family} 

The discussion in this section is very similar to the corresponding
discussion for $K3$ surfaces given in section 3 of \cite{Br-Le}. We show
that there is a unique family $T$ (up to deformation) corresponding to the
twistor family. We define $N_{g,n}^{FLS}$ and $N_{g,n}$ using a suitable
set of Gromov-Witten invariants for the family $T$. We show that the
invariants solve the enumerative problems we are interested in. We use this
framework to prove that the invariants $N^{FLS}_{g,n}$ and $N_{g,n}$ are
universal numbers independent of the choice of the Abelian surface $X$ and
the linear system $|L|$.

Let $(X,\omega )$ be an Abelian surface and let $g$ be the unique
hyperk\"ahler metric given by Yau's theorem \cite{Yau}. Define $T_{g}$ to
be the family of K\"ahler structures given by the unit sphere in
$\mathcal{H}^{2}_{+,g}$, the space of self-dual, harmonic 2-forms.

\begin{prop}\label{prop: twistor family is unique}
For any two hyperk\"ahler metrics $g$ and $g'$, the corresponding twistor
families $T_{g}$ and $T_{g'}$ are deformation equivalent. There is
therefore a well-defined deformation class which we denote by $T$.
\end{prop}
\pf The moduli space of complex structures on the 4-torus is connected and
the space of hyperk\"ahler structures for a fixed complex torus is
contractible (it is the K\"ahler cone). Therefore the space parametrizing
hyperk\"ahler 4-tori is connected (in fact it is just the space of flat
metrics). We can thus find a path of hyperk\"ahler metrics connecting $g$
to $g'$ and by associating the twistor family to each metric, we obtain a
continuous deformation of $T_{g}$ to $T_{g'}$. \qed 

An important observation concerning the twistor family is the following
corollary.
\begin{cor}
Let $f:X\to X$ be an orientation preserving diffeomorphism, then $f^{*}(T)$
is deformation equivalent to $T$. 
\end{cor}
\pf Let $T_{g}$ be the twistor family for a hyperk\"ahler metric $g$. Then
$f^{*}(T_{g})=T_{f^{*}(g)}$ where $f^{*}(g)$ is the pullback metric which
is also hyperk\"ahler. \qed 

The corollary has the consequence that for any orientation preserving
diffeomorphism $f$ the Gromov-Witten invariants for the twistor family
satisfy 
$$
\Phi ^{X,T}_{g,C}(\alpha _{1},\ldots,\alpha _{l})=\Phi
^{X,T}_{g,f_{*}(C)}(f_{*}\alpha _{1},\ldots,f_{*}\alpha _{l}).
$$

There is an orientation diffeomorphism of the 4-torus for every element of
$Sl_{4}(\znums )$ given by the descent of the linear action on the
universal cover $\rnums ^{4}$. It follows from the elementary divisor
theorem that for any 
two classes $C $ and $C' $ in $H_{2}(X,\znums )$ with the same divisibility
and square there is a linear diffeomorphism $f$ such that $f_{*}(C)=C'$
(\cite{abelian} pg. 47). This means that there is a lot of symmetry among
the Gromov-Witten invariants for the twistor family so that (for primitive
classes) they can be  encompassed by $N_{g,n} $ and $N_{g,n}^{FLS}$ and
the other invariants which we define below.

\begin{definition}\label{def: definition of Ngn and NngFLS}
Let $\gamma _{1},\ldots,\gamma _{4}$ be loops in $X$ representing an
oriented basis for 
$H_{1}(X,\znums )$ and let $C\in H_{2}(X,\znums )$ be a primitive class
with $C\cdot C=2g-2+2n$. We define:
\begin{eqnarray*}
N_{g,n}&=&\Phi ^{(X,T)}_{g,C}(\point ^{g})\\
N_{g,n}^{FLS}&=&\Phi ^{(X,T)}_{g,C}(\gamma _{1},\gamma _{2},\gamma
_{3},\gamma _{4},\point ^{g-2}).
\end{eqnarray*}
Now let $C=[\gamma _{1}]\wedge  [\gamma _{2}]+(n+g-1)[\gamma _{3}]\wedge
[\gamma _{4}] $ and define\footnote{Here we use the fact that on any torus
there is a 
natural identification $$\Lambda ^{2}H_{1}(X;\znums )\cong H_{2}(X,\znums
).$$ }   
$$
N^{ij}_{g,n}=\Phi ^{(X,T)}_{g,C}(\gamma _{i},\gamma _{j},\point ^{g-1})
$$
for $i<j$.
\end{definition}

The invariants $N_{g,n}$, $N_{g,n}^{FLS}$, and $N_{g,n}^{ij}$ encompass all
possible Gromov-Witten invariants for the twistor family and primitive
homology classes. Since $X$ has real dimension four, the only non-trivial
constraints come 
from the point class and $H_{1}(X;\znums )$. Since the point class is
invariant under orientation preserving diffeomorphisms, $\Phi
^{(X,T)}_{g,C}(\point ^{g})$ only depends on the square (and divisibility)
of $C$. Similarly, since $\Phi $ is skew-symmetric in the classes $\gamma
_{i}$, the only possibility with four
 $\gamma $-constraints are the invariants $\Phi _{g,C}^{(X,T)}(\gamma
_{1},\gamma _{2},\gamma _{3},\gamma _{4},\point ^{g-2})$ which also only
depends on the square (and divisibility) of $C$.
\footnote{This invariant only depends on the square and divisibility of $C$
since for any orientation preserving diffeomorphism $f$ we have:
\begin{eqnarray*}
\Phi ^{(X,T)}_{g,C}(\gamma _{1},\ldots,\gamma _{4},\point ^{g-2})&=&
\Phi _{g,f_{*}(C)}^{(X,T)}(f_{*}\gamma _{1},\ldots,f_{*}\gamma _{4},\point
^{g-2})\\
&=&\det (f_{*}:H_{1}\to H_{1})\Phi _{g,f_{*}(C)}^{(X,T)}(\gamma
_{1},\ldots,\gamma _{4},\point ^{g-2})\\
&=&\Phi _{g,f_{*}(C)}^{(X,T)}(\gamma
_{1},\ldots,\gamma _{4},\point ^{g-2}).
\end{eqnarray*}
} 

For dimensional reasons, the invariants with one or three $\gamma $'s
 are zero. The invariants $N_{g,n}^{ij}$ encompass the remaining invariants
since for any primitive $C$ with $C^{2}=2g-2+2n$ we can first move $C$ to
$[\gamma _{1}]\wedge  [\gamma _{2}]+(n+g-1)[\gamma _{3}]\wedge
[\gamma _{4}] $ by an orientation preserving diffeomorphism. 

We now wish to show that $N_{g,n}$, $N_{g,n}^{FLS}$, and $N_{g,n}^{ij}$
enumerate holomorphic curves as we want.

\begin{lemma}\label{lem: twistor family has only one member supposrting
curves in C}
Suppose that $X$ is an Abelian  surface with a hyperk\"ahler metric $g$ and
suppose that $C\subset X$ is a holomorphic curve. Then the only K\"ahler
structure in $T_{g}$ that has a holomorphic curve in the class $[C]\in
H_{2}(X;\znums )$ is the original K\"ahler structure for which $C$ is holomorphic.
\end{lemma}
\pf A necessary condition for the class $[C]$ to contain holomorphic curves
is that $[C]^{\vee }\in H^{1,1}(X,\rnums )$ and $[C]$ pairs positively with
the K\"ahler form. Since the class orthogonal to $\mathcal{H}^{2}_{+,g}$
are always of type $(1,1)$ we just need to see when the projection of
$[C]^{\vee }$
to $\mathcal{H}^{2}_{+,g}$ is type $(1,1)$.  

We may assume $[C]^{2}$ is non-negative since $X$ has no rational curves
and so the projection of $[C]^{\vee }$ onto $\mathcal{H}^{2}_{+,g}$ is
non-zero. By 
definition, the complex structure associated to a form $\omega \in
T_{g}\subset \mathcal{H}^{2}_{+,g}$ defines the orthogonal splitting
$\mathcal{H}^{2}_{+,g}\cong \omega \rnums \oplus (H^{2,0}\oplus
H^{0,2})_{\rnums }$. Therefore, the only $\omega \in T_{g}$ for which
$[C]^{\vee  }$
is type $(1,1)$ and pairs positively with $\omega $ is when $\omega $ is a
positive multiple of the projection of $[C]^{\vee }$ to
$\mathcal{H}^{2}_{+,g}$. 
This is unique and so must be the original K\"ahler structure for which
$C\subset X$ is holomorphic. \qed 

In general, Gromov-Witten type invariants may give a different count from
the corresponding enumerative problem because Gromov-Witten invariants
count maps instead of curves and the image of a map may have geometric
genus smaller than its domain. We show that in the case at hand this does
not happen.
\begin{lemma}
Suppose $X$ is generic among those Abelian surfaces admitting a curve in
the class of $C$. Then the invariants $N_{g,n}$, $N_{g,n}^{FLS}$, and
$N_{g,n}^{ij}$ count only maps whose image has geometric genus $g$ and they
are counted with positive integral multiplicity.
\end{lemma}
\pf The assumption guarantees that $C$ generates $\pic (X)/\pic ^{0}(X)$ so
that all the curves in the class $[C]$ are reduced and irreducible. Then a
dimension count shows that maps with contracting components of genus one or
greater are not counted. The details are the same as in the proof of
Theorem 3.5 in \cite{Br-Le}. \qed

We summarize the results of this section in the following theorem.

\begin{theorem}\label{thm: summary of properties of invariants}
Let $C\subset X$ be a holomorphic curve $C$ in an Abelians surface $X$ with
$C^{2}=2g-2+2n$. Then the invariants $N_{g,n}$, $N_{g,n}^{FLS}$, and
$N^{ij}_{g,n}$ defined in Definition \ref{def: definition of Ngn and
NngFLS} count the number of genus $g$ curves in the class $[C]$
that:
\begin{enumerate} 
\item pass through $g$ generic points ($N_{g,n}$);
\item pass through $g-2$ generic points and are confined
to the linear system $|C|$ ($N_{g,n}^{FLS}$). 
\item pass though $g-1$ generic points and are confined to lie in a certain
1-dimensional family of linear systems determined by $ij\in
\{12,13,14,23,24,34 \}$ ($N_{g,n}^{ij}$).  
\end{enumerate}
Furthermore, the invariants do not depend on the particular Abelian surface
$X$ or the curve $C$ and the invariants count all irreducible, reduced
curves with positive multiplicity that is one for nodal curves. 
\end{theorem}
\section{Computing the invariants}\label{sec: computing Ng(n)}

In this section, we compute the invariants $N_{g,n}$, $N_{g,n}^{FLS}$ and
$N_{g,n}^{ij}$  by using a particular choice for $X$ and $C$. Namely, we
will choose $X$ to be a product of elliptic curves and $C$ to be a section
together with a multiple of the fiber. The computations 
give our main theorem as well as the generalization to the $N_{g,n}^{ij}$
case.  For convienence we state the results below:

\begin{theorem}\label{thm: restatement of main thm including Nij}
The invariants  $N_{g,n}$, $N_{g,n}^{FLS}$ and
$N_{g,n}^{ij}$  (defined in Definiton \ref{def:
definition of Ngn and NngFLS}) are given by the following generating
functions: 
\begin{eqnarray}
 \label{eqn: formula for Ngn}
\sum_{n=0}^\infty N_{g,n}q^{n+g-1} &=&g\left( DG_2\right) ^{g-1}\\
 \label{eqn: formula for
NFLS} 
\sum_{n=0}^{\infty }N_{g,n}^{FLS}q^{n+g-1}
&=&(DG_{2})^{g-2}D^{2}G_{2}\\
&=&(g-1)^{-1}D((DG_{2})^{g-1})\nonumber\\
 \label{eqn: formula for N12} 
\sum_{n=0}^\infty N_{g,n}^{12}q^{n+g-1}&=&D\left((DG_{2})^{g-1} \right)\\
\label{eqn: formula for N34} 
\sum_{n=0}^\infty N_{g,n}^{34}q^{n+g-1}&=&\left( DG_2\right) ^{g-1}\\
 \label{eqn: N13,N14,N23 are 0} 
\sum_{n=0}^{\infty }N_{g,n}^{13}q^{n+g-1}&=&0\\
\sum_{n=0}^{\infty
}N_{g,n}^{14}q^{n+g-1}&=&0\nonumber\\
\sum_{n=0}^{\infty
}N_{g,n}^{23}q^{n+g-1}&=&0\nonumber\\
\sum_{n=0}^{\infty
}N_{g,n}^{24}q^{n+g-1}&=&0.\nonumber
\end{eqnarray}
where $D=q\frac{d}{dq}$ and $G_{2}=-\frac{1}{24}+\sum_{k=1}^{\infty
}(\sum_{d|k}d)q^{k}$. 
\end{theorem}

From Theorem \ref{thm: summary of properties of invariants} we are free to
compute the invariants for any choice of $X$ and $C$.  Let $X$ be the
product of two generic elliptic 
curves $S\times F$ and let $C$ be the primitive homology class $S+\left(
g+n-1\right) F\in H_2\left( X,
\mathbf{Z}\right)$.  We write $X=S^1\times S^1\times S^1\times S^1$ by
choosing 
a diffeomorphism between $S^1\times S^1$ with $S$ and also one between $%
S^1\times S^1$ with $F$ . Next we choose representatives for the loops and
points on $X.$ We consider following four loops generating $H_1\left( X,%
\mathbf{Z}\right) :$

\begin{equation*}
\begin{array}{ll}
\gamma _1:S^1\rightarrow X, & \gamma _1\left( e^{it}\right) =\left(
e^{it},b_1,c_1,d_1\right) , \\
\gamma _2:S^1\rightarrow X, & \gamma _2\left( e^{it}\right) =\left(
a_2,e^{it},c_2,d_2\right) , \\
\gamma _3:S^1\rightarrow X, & \gamma _3\left( e^{it}\right) =\left(
a_3,b_3,e^{it},d_3\right) , \\
\gamma _4:S^1\rightarrow X, & \gamma _4\left( e^{it}\right) =\left(
a_4,b_4,c_4,e^{it}\right) .
\end{array}
\end{equation*}
where the $a_{i}$'s, $b_{i}$'s, $c_{i}$'s, and $d_{i}$'s  are distinct
points on $S^1.$ We also choose $g$ 
generic points $p_1=\left( s_1,f_1\right) ,...,p_g=\left( s_g,f_g\right) $
on $X=S\times F.$ For any $f\in F$ and $s\in S$ we call $S_f=S\times \left\{
f\right\} \subset X$ a section curve and $F_s=\left\{ s\right\} \times F$ a
fiber curve of $X.$

{\sc Proof of Equation \ref{eqn: formula for Ngn}:}

We suppose that $\phi :D\rightarrow X$ is a stable map from a genus $g$
curve $D$ with $g$ marked points to $X$ representing the class $C=S+\left(
g+n-1\right) F$ and sending corresponding marked points to the $p_{i}$'s.
First we observe that the image of $\phi $ consists of one section curve 
and some fiber curves. This is because the projection of any irreducible
component of $\im (\phi ) $ to $S$ is of degree zero except for
one which has degree one. On the other hand, the projection of the degree
one component to $F$ must have degree
zero because there is no nontrivial morphism between two generic elliptic
curves $S$ and $F$.  Therefore $\im (\phi )$ consists of a single
section curve  and a number of fibers.

In order for $\phi $ to pass through the $g$ generic points
$p_{1},\ldots,p_{g}$ 
on $X,$ we need at least $\left( g-1\right) $ fiber curves in the image of $%
\phi .$ On the other hand the geometric genus of $D\ $is $g$ and covering
each section or fiber curve will take up at least one genus of $D$. 
Therefore $\im \left( \phi \right) $ consists of exactly $g-1$ fiber curves
(possibly with multiplicity) and one section curve.

In fact, $D$ has to have precisely $g$ irreducible components
$D_1,\ldots,D_g$ and 
each component is a genus one curve and contains one marked point. The
restriction of $\phi $ to each $D_i$ is either (i) a covering of some fiber
curve containing one of the $p_{i}$'s or (ii) an isomorphism to a
section curve containing one of the $p_{i}$'s.

We can assume that $\phi $ restricted to $D_g$ is an isomorphism onto one of
 the section curves containing some $p_i=\left( s_i,f_i\right) .$ There are
 $g$ 
choices of such $p_{i}$'s. Without loss of generality we assume $\phi
\left( D_g\right) $ contains $p_g$ or equivalently $\phi \left( D_g\right)
=S_{f_g}.$ In this case the marked point on $D_g$ must be the unique point
 $\phi ^{-1}\left( p_g\right)$. 

For $i<g,$ we label the component of $D$ covering $F_{s_i}$ as $D_i.$ Then
$\phi :D_i\rightarrow F_{s_i}$ is an unbranched cover by the Hurwitz
theorem 
because both $D_i$ and $F_{s_i}$ are genus one curves. We denote
$k_i=\operatorname{deg}\left( \phi :D_i\rightarrow F_{s_i}\right) >0$ then
we have 
$\sum_{i=1}^{g-1}k_i=g+n-1$ since $\phi $ represents the homology class of
$S+\left( g+n-1\right) F$. The number of elliptic curves that admit a
degree $k$ homomorphism to  a fixed
elliptic curve is classically known to be $\sum_{d|k}d$. We
fix the origin of the elliptic curve $D_i$ (and $F_{s_i}$) to be its
intersection point with $D_g$ (and $S_{f_g}$ respectively), so that $\phi
:D_{i}\to F_{s_{i}}$ is a homomorphism and thus the
number of choices of $\left( \phi :D_i\rightarrow F_{s_i}\right) $ is given
by $\sum_{d|k_i}d$.  Since the marked point on $D_i$ could be any one of the
$k_i$ points in $\phi ^{-1}\left( p_i\right)$, there are a total of
$k_{i}\sum_{d|k_{i}}d$ choices for each marked curve $D_{i}$. 

We denote the $(g-1)$-tuple $k_{1},\ldots,k_{g-1}$ by $\mathbf{k}$ and we
write $|\mathbf{k}|$ for $\sum_{i}k_{i}$. 
From the preceding discussion, the number of stable maps $\phi $ of
geometric genus $g$ and $g 
$ marked points to $X$ in the class of $S+\left( g+n-1\right) F$ is given by
\begin{equation*}
g\sum_{\mathbf{k}:|\mathbf{k}|=n+g-1} \prod_{i=1}^{g-1}k_i\left(
\sum_{d|k_i}d\right).
\end{equation*}
It is 
not difficult to see that each such stable map $\phi $ contributes one to
the family Gromov-Witten invariant (c.f. \cite{Br-Le}).
The formula is then proved by summing over $n$ and rearranging:
\begin{eqnarray*}
\sum_{n=0}^\infty N_{g,n}q^{n+g-1} &=&\sum_{n=0}^\infty \left(
g\sum_{|\mathbf{k}|=n+g-1}\prod_{i=1}^{g-1}k_i\left(
\sum_{d|k_i}d\right) \right) q^{n+g-1} \\
&=&g\sum_{n=0}^\infty \left(
\sum_{|\mathbf{k}|=n+g-1}\prod_{i=1}^{g-1}\left(
k_i\sum_{d|k_i}dq^{k_{i}}\right)\right)  \\
&=&g\left( \sum_{k=1}^\infty k\sum_{d|k}dq^k\right) ^{g-1}\\
&=&g(DG_{2})^{g-1}.
\end{eqnarray*}\qed

{\sc Proof of Equations  \ref{eqn: formula for N12},  \ref{eqn: formula for
N34}, and \ref{eqn: N13,N14,N23 are 0}:} 

For these computations we count the number of stable maps $\phi
:D\rightarrow X$ of genus $g$ with $g+1$ marked points that represent the
homology class $C=S+\left( g+n-1\right) F$. The maps are constrained by
requiring that the first $g-1$ marked points are mapped to the points
$p_{1},\ldots,p_{g-1}$ and the image of the remaining two points must lie
on $\gamma _{i}$ and $\gamma _{j}$ respectively. We argue as before that
$D$ has $g$ irreducible components $D_1,...,D_g$  where
each $D_i$ is a genus one curve and it contains one marked point except $D_g$
which contains two marked points. Moreover, the image of the two marked
points on $D_g$ must lie on the two loops $\gamma _{i}$ and $\gamma _{j}$.
It is easy 
to check that when the points $p_{i},a_{i},\ldots,d_{i}$  are in general
positions, then no single section or fiber curve can pass through both one
of the loops $ \gamma_{i}$ and one of the points $p_{i}$. We can also verify
directly the following lemma.

\begin{lemma}
No fiber or section curve can pass through two different
$\gamma_{i}$'s unless they are $\left( \gamma _1,\gamma _2\right) $ or
$\left( \gamma _3,\gamma _4\right)$.

No section curve can pass through both $\gamma _1$ and
$\gamma _2.$ The only fiber curve passing through $\gamma _{1}$ and $\gamma
_{2}$  is $F_{12}=F_{\left( a_2,b_1\right) }$.
Moreover, $F_{12}\cap \gamma _1 =\left(
a_2,b_1,c_1,d_1\right) $ and $F_{12}\cap \gamma _2 =\left(
a_2,b_1,c_2,d_2\right) .$

Similarly,  no fiber curve can pass through both $\gamma _3$
and $\gamma _4.$ The 
only section curve passing through $\gamma _{3}$ and $\gamma _{4}$  is
$S_{34}=S_{\left( c_4,d_3\right) }$.  Moreover, $S_{34}\cap  \gamma _3
=\left( a_3,b_3,c_4,d_3\right) $ and $S_{34}\cap \gamma _4 =\left(
a_4,b_4,c_4,d_3\right)$. 
\end{lemma}

From the lemma the moduli space of stable maps is empty in all cases except $%
\left( \gamma _i,\gamma _j\right) =\left( \gamma _1,\gamma _2\right) $ or $%
\left( \gamma _3,\gamma _4\right) $ up to permutations. Therefore the
corresponding family Gromov-Witten invariants vanish (proving Equation
\ref{eqn: N13,N14,N23 are 0}) 
\begin{equation*}
N_{g,n}^{13}=N_{g,n}^{14}=N_{g,n}^{23}=N_{g,n}^{24}=0.
\end{equation*}

Next we compute $N_{g,n}^{34}.$ Let $\phi $ be any stable map in the
corresponding moduli space, we have $\phi \left( D_g\right) =S_{34}$ by the
above lemma. Moreover the restriction of $\phi $ to $D_g$ is an isomorphism
and $\phi $ must send the two marked points on $D_g$ to $\left(
a_3,b_3,c_4,d_3\right) $ and $\left( a_4,b_4,c_4,d_3\right) .$ It is not
difficult to check that the orientation of the moduli space is compatible
with the one induced from $\left( \gamma _3,\gamma _4\right) .$ The
calculation for the contribution from the other $D_i$'s are identical
with the earlier computation and we obtain 
$$
\sum_{n=0}^\infty N_{g,n}^{34}q^{n+g-1}=\left( \sum_{k=1}^\infty
k\sum_{d|k}dq^k\right) ^{g-1}=(DG_{2})^{g-1}
$$
thus proving Equation \ref{eqn: formula for N34}.

Now we compute $N_{g,n}^{12}$. From the lemma again we have $\phi \left(
D_g\right) =F_{12}$ for any stable map $\phi $ in the corresponding moduli
space. Let us denote the degree of $\phi $ restricted to $D_g$ by $k_0.$
Then the number of possible $\left( \phi :D_g\rightarrow F_{12}\right) $ is
given by $\sum_{d|k_0}d$ as before. The image of the two marked points on $%
D_g$ are $\left( a_2,b_1,c_1,d_1\right) $ and $\left( a_2,b_1,c_2,d_2\right)
.$ Since $\phi $ is an unbranched covering map, there are $k_0$ choices for
the location of each marked point. There are then a total of
$k_0^2\sum_{d|k_0}d$ different choices associated the component $D_g$. 

There are $g-1$ choices for which curve $D_{1},\ldots,D_{g-1}$ is mapped to
a section. After making such a choice we may assume without loss of
generality that $D_{g-1}$ is mapped to the section curve passing through
$p_{g-1}$  (since we can relabel the points). Then since $\phi $ is an
isomorphism on $D_{g-1}$, the location of the 
marked point on it is determined. The analysis for the other fibers is
the same as before and the issue of the orientation is the same as above.
In conclusion we have
\begin{equation*}
N_{g,n}^{12}=\left( g-1\right) \sum_{\mathbf{k}:|\mathbf{k}|=n+g-1}  \left(
k_0^2\sum_{d|k_0}d\prod_{i=1}^{g-2}k_i\sum_{d|k_i}d\right) ,
\end{equation*}
where the summation is over $(g-1)$-tuples
$\mathbf{k}=\{k_{0},\ldots,k_{g-2}  \}$ with $k_i>0$ and
$|\mathbf{k}|=\sum_{i=0}^{g-2}k_i=g+n-1$. Summing over $n$, we prove
Equation \ref{eqn: formula for N12}:

\begin{eqnarray*}
\sum_{n=0}^\infty N_{g,n}^{12}q^{n+g-1} &=&\left( g-1\right) \left(
\sum_{k=1}^\infty k^2\sum_{d|k}dq^k\right) \left( \sum_{k=1}^\infty
k\sum_{d|k}dq^k\right) ^{g-2} \\
&=&\left( g-1\right) \left( D^2G_2\right) \left( DG_2\right) ^{g-2}\\
&=&D\left((DG_{2})^{g-1} \right).
\end{eqnarray*}\qed

{\sc Proof of Equation \ref{eqn: formula for
NFLS}:} 

Suppose that $\phi :D\rightarrow X$ is a stable map of genus $g$ with $g+2$
marked points and represents the homology class $C=S+\left( g+n-1\right) F$.
We argue as before that $D$ has $g$ irreducible components $D_1,...,D_g$.
The curves $D_{1},\ldots,D_{g-2}$ each have one marked point mapping to
$p_{1},\ldots,p_{g-2}$ respectively. The curves $D_{g-1}$
and $D_g$ contain the remaining four marked points which must lie on the
four loops $\gamma 
_1,...,\gamma _4$. The image of each $D_i$ must either be a section curve or
a fiber curve. By the previous lemma the only fiber curve that intersects
two $\gamma_{i}$'s must be $F_{12}$ and it intersects $\gamma _1$ and
$\gamma _2$ 
at $\left( a_2,b_1,c_1,d_1\right) $ and $\left( a_2,b_1,c_2,d_2\right) $
respectively.
Similarly, the only section that intersects two $\gamma$'s is $S_{34} $ and
it must intersect $\gamma _3$ and $\gamma _4$ at $\left( 
a_3,b_3,c_4,d_3\right) $ and $\left( a_4,b_4,c_4,d_3\right) $
respectively. 
Hence the image of $D_{g-1}$ and $D_g$ must be $F_{12}$ and $S_{34}$
respectively. For $i\leq g-2,$ we have $\phi \left( D_i\right) =F_{s_i}.$
Let 
the degree of $\phi $ on $D_i$ be $k_i$ where $1\leq i\leq g-1.$ Then $\sum
k_i=g+n-1.$ The number of choices of $\left( \phi :D_i\rightarrow
F_{s_i}\right) $ is again given by $\sum_{d|k_i}d.$ For $i<g-1$ there are $%
k_i$ choices of the location of the marked point on $D_i.$ For $i=g-1,$
there are two marked points on $D_{g-1}$ and each has $k_{g-1}$ possible
locations. Again it is not difficult to check that such maps each
contribute $+1$ to the family Gromov-Witten invariants. In conclusion we have
\begin{equation*}
N_{g,n}^{FLS}=\sum_{\mathbf{k}:|\mathbf{k}|=n+g-1} \left(
k_{g-1}^2\sum_{d|k_{g-1}}d\right) 
\prod_{i=1}^{g-2}\left( k_i\sum_{d|k_i}d\right)
\end{equation*}
and so summing over $n$ we have
\begin{eqnarray*}
\sum_{n=0}^\infty N_{g,n}^{FLS}q^{n+g-1} &=&\left( \sum_{k=1}^\infty
k\sum_{d|k}dq^k\right) ^{g-2}\left( \sum_{k=1}^\infty
k^2\sum_{d|k}dq^k\right)  \\
&=&\left( DG_2\right) ^{g-1}\left( D^2G_2\right) 
\end{eqnarray*}
proving Equation \ref{eqn: formula for
NFLS}.\qed

\appendix
\section{Proof of Theorem \ref{thm: pullback of class on pic is
ft*ev*}}\label{appendix}  

First we determine a more explicit formulation of the map
$$
\Psi _{\Sigma _{0}}:\M _{g,C}\to \pic ^{0}(X)=H^{1}(X,\rnums
)/H^{1}(X,\znums ).
$$
Throughout this section we will abuse notation and refer to a stable map
$f:\Sigma \to X$ by its image so that when we say $\Sigma \in \M_{g,C}$ we
mean the curve $\Sigma \subset X$ given by the image of the map $f:\Sigma
\to X$. By definition, $\Psi_{\Sigma _{0}} (\Sigma )=\mathcal{O}(\Sigma
-\Sigma _{0})$. We wish to get a 1-form representative for
$\mathcal{O}(\Sigma -\Sigma _{0})$. First, choose a $C^{\infty }$
trivialization of the bundle $\mathcal{O}(\Sigma -\Sigma _{0})$. Let
$s_{0}$ and $s$ be defining sections of $\mathcal{O}(\Sigma _{0})$ and
$\mathcal{O}(\Sigma )$. Then via the trivialization, $h=s_{0}/s$ is a
non-vanishing $C^{\infty }$ function on $X^{o}=X-\{\Sigma _{0}\cup \Sigma
\}$. Define an integrable real  1-form $a_{h}$ by
$$
a_{h}=\frac{1}{2\pi i}\left(\frac{\dbar h}{h}-\frac{\del
\overline{h}}{\overline{h}} \right).
$$
The form $a_{h}$ defines  a real 1-current and let $a_{\Sigma _{0}}^{\Sigma
}$ be the 1-form 
representing the harmonic projection of $a_{h}$. We wish to define a map to
$H^{1}(X,\rnums )/H^{1}(X,\znums )$ using $a_{\Sigma _{0}}^{\Sigma }$ and
we need the following lemma:
\begin{lemma}
The harmonic 1-form  $a_{\Sigma _{0}}^{\Sigma }$ is independent of the
choices of $s_{0}$, $s$, and the trivialization up to translation by an
integral harmonic 1-form.
\end{lemma}
\pf Let $h'=s_{0}'/s'$ be defined using a possibly different
trivialization. Then $h'=gh $ for some $C^{\infty }$ map $g:X\to \cnums
^{*}$ and so 
$$
a_{h}-a_{h'} = \frac{1}{2\pi i}\left(\frac{\dbar g}{g}-\frac{\del
\overline{g}}{\overline{g}} \right) = a_{g}
$$
is a smooth 1-form. We need to show that the harmonic projection of $a_{g}$
is integral. Write $g=e^{t}\theta $ where $t:X\to \rnums $ and $\theta
:X\to S^{1}\subset \cnums $ are $C^{\infty }$. A direct computation shows that 
$$
a_{g}=d^{c}t+\frac{1}{2\pi i}\theta ^{-1}d\theta =d^{*}(t\omega
)+\frac{1}{2\pi i}\theta ^{-1}d\theta 
$$
where the last equality uses the K\"ahler identities and so we see that
the harmonic projection of $a_{g}$ is the integral class $\frac{1}{2\pi
i}\theta ^{-1}d\theta $.\qed 

The lemma shows that 
$$
\Sigma \mapsto a_{\Sigma _{0}}^{\Sigma }
$$
determines a well-defined map
$$
\M _{g,C}\to H^{1}(X,\rnums )/H^{1}(X,\znums ).
$$
This is the same map as $\Psi _{\Sigma _{0}}$ since
\begin{enumerate}
\item if $\Sigma $ is linearly equivalent to $\Sigma _{0}$ then $a_{\Sigma
_{0}}^{\Sigma }\equiv 0$ and
\item the sum of divisors corresponds to the sum of forms, \ie $a_{2\Sigma
_{0}}^{\Sigma +\Sigma '}=a_{\Sigma _{0}}^{\Sigma 
}+a_{\Sigma _{0}}^{\Sigma '}$. 
\end{enumerate}
The first property holds since if $\Sigma \sim \Sigma _{0}$ then $h$ can be
chosen to be a meromorphic function and so $a_{h}\equiv 0$. The second
property follows from $a_{h h'}=a_{h}+a_{h'}$.

To prove Theorem \ref{thm: pullback of class on pic is ft*ev*} we need  to
show that $\Psi ^{*}_{\Sigma _{0}}(\til{\gamma 
} )=ft_{*}ev^{*}([\gamma  ]^{\vee })$. Recall that $\gamma $ is a loop in
$X$, $[\gamma ]^{\vee}$ is the Poincare dual of the 1-cycle $[\gamma ]$,
and $\til{\gamma  }\in H^{1}(\pic ^{0}(X),\znums )$ is the natural class
associated to $[\gamma ]$ arising from the identification $\pic
^{0}(X)\cong H^{1}(X,\rnums )/H^{1}(X,\znums )$. 

We use the general correspondence 
between elements of $H^{1}(M,\znums )$ and homotopy classes of maps $M\to
S^{1}$. One can 
get a circle valued function on $M$ from a class $\phi \in H^{1}(M;\znums
)$ by choosing a base point $x_{0}\in M$ and defining $f_{\phi }:M\to
\rnums /\znums $ by
$$
f_{\phi }(x)=\int_{\Gamma _{x_{0}}^{x}}\phi \mod \znums 
$$
where $\Gamma _{x_{0}}^{x}$ is a path from $x_{0}$ to $x$. Since $\phi $ is
an integral class, $f_{\phi }$ does not depend on the choice of the path
(mod $\znums $).

Using $\Sigma _{0}$ as the base point for $\M _{g,C}$, the class
$ft_{*}ev^{*}([\gamma  ]^{\vee })\in H^{1}(\M _{g,C},\znums )$ is given by the
$S^{1}$-valued map
\begin{eqnarray*}
\Sigma &\mapsto &\int_{\Gamma _{\Sigma _{0}}^{\Sigma
}}ft_{*}ev^{*}([\gamma  ]^{\vee }) \mod \znums \\
&=&\int_{ev(ft^{-1}(\Gamma _{\Sigma _{0}}^{\Sigma }))}[\gamma  ]^{\vee } \mod
\znums \\
&=&\int_{W_{\Sigma _{0}}^{\Sigma }}[\gamma  ]^{\vee }\mod \znums 
\end{eqnarray*}
where $W_{\Sigma _{0}}^{\Sigma }$ is the 3-chain in $X$ swept out by the
curves in the path $\Gamma _{\Sigma _{0}}^{\Sigma }$. Note that $\del
W_{\Sigma _{0}}^{\Sigma }=\Sigma _{0}-\Sigma $ and that the map
$$
\Sigma \mapsto \int_{W}[\gamma  ]^{\vee }\mod \znums 
$$
is the same for any 3-chain $W$ such that $\del W=\Sigma _{0}-\Sigma $
(since the difference $W_{\Sigma _{0}}^{\Sigma }-W$ is a 3-cycle and
$[\gamma  ]^{\vee }$ is an integral class).

On the other hand, the class $\til{\gamma } \in H^{1}(\pic ^{0},\znums
)$ is by definition given by the $S^{1}$-valued function on $\pic ^{0}$
defined by 
$$
[a]\mapsto \int_{X}a\wedge [\gamma  ]^{\vee } \mod \znums 
$$
where $a\in H^{1}(X,\rnums )$ and $[a]$ is the corresponding equivalence
class in 
$$H^{1}(X,\rnums )/H^{1}(X,\znums ).$$ 
The class $\Phi _{\Sigma
_{0}}^{*}(\til{\gamma } )$ is therefore represented by the
$S^{1}$-valued function
$$
\Sigma \mapsto \int_{X}a_{\Sigma _{0}}^{\Sigma }\wedge [\gamma  ]^{\vee } \mod
\znums 
$$
and to prove theorem \ref{thm: pullback of class on pic is ft*ev*} then we need
to show that 
$$
\int_{X}a_{\Sigma _{0}}^{\Sigma }\wedge [\gamma  ]^{\vee }=\int_{W_{\Sigma
_{0}}^{\Sigma }} [\gamma  ]^{\vee }\mod \znums .
$$
Recall that $a_{\Sigma _{0}}^{\Sigma }$ is the harmonic projection of
$a_{h}$. Let $\zeta  $ be the harmonic representative for $[\gamma
 ]^{\vee }$. Then 
$$
\int_{X}a_{\Sigma _{0}}^{\Sigma }\wedge [\gamma  ]^{\vee }=\int_{X}a_{h}\wedge
\zeta  
$$ 
so we need to show that
\begin{equation}\label{eqn: int over Xo is same as in over W mod Z}
\int_{X^{o}}\frac{1}{2\pi i}\left(\frac{\dbar h}{h}-\frac{\del
\overline{h}}{\overline{h}} \right)\wedge \zeta  =\int_{W_{\Sigma
_{0}}^{\Sigma }}\zeta   \mod \znums .
\end{equation}
Recall that $h$ is a $\cnums ^{*}$-valued function on $X^{o}=X-\{\Sigma
_{0}\cup \Sigma  \}$ and so writing $h=e^{t}\theta $ for smooth functions
$t:X^{o}\to \rnums $ and $\theta :X^{o}\to S^{1}$ we can rewrite the left
hand side of Equation \ref{eqn: int over Xo is same as in over W mod Z} as 
$$
LHS=\int_{X^{o}}d^{*}(t\omega )\wedge \zeta  +\int_{X^{o}}\frac{1}{2\pi
i}\theta ^{-1}d\theta \wedge \zeta  .
$$
We can do the first of these integrals by first integrating along the
fibers of $t:X^{o}\to \rnums $:
\begin{eqnarray*}
t_{*}(d^{*}(t\omega )\wedge \zeta  )&=&-t_{*}(*d(*(t\omega ))\wedge \zeta
 ) \\
&=&t_{*}(d(t\omega )\wedge *\zeta  )\\
&=&t_{*}(dt\wedge \omega \wedge *\zeta  )\\
&=&dt\wedge t_{*}(\omega \wedge *\zeta  )
\end{eqnarray*}
but $t_{*}(\omega \wedge *\zeta  )=0$ because $\omega \wedge *\zeta
 $ is a closed form defined on {\em all} of $X$ and $t^{-1}(\point )$ is
a boundary 3-chain in $X$.

We perform the remaining integral by integrating first along the fibers of
$\theta :X^{o}\to S^{1}$. Since $\frac{1}{2\pi i}\theta ^{-1}d\theta=\theta
^{*}(dvol_{S^{1}})$, we can write the remaining integral as 
$$
\int_{X^{o}}\theta ^{*}(dvol_{S^{1}})\wedge \zeta  =\int_{S^{1}}\theta
_{*}(\zeta  ) dvol_{S^{1}}.
$$
Now $\theta _{*}(\zeta  )=\int_{\theta ^{-1}(c)}\zeta  $ is
independent of $c$ (mod $\znums $) since $\theta ^{-1}(c)$ is a 3-chain in
$X$ with boundary $\Sigma _{0}-\Sigma $ and so (mod $\znums $) we have
\begin{eqnarray*}
\text{LHS of Eqn \ref{eqn: int over Xo is same as in over W mod Z} }
&=&\left(\int_{W}\zeta   \right)\int_{S^{1}}dvol_{S^{1}}\\
&=&\int_{W}\zeta  \\
&=&\text{RHS of Eqn \ref{eqn: int over Xo is same as in over W mod Z} } 
\end{eqnarray*}
which proves the theorem. \qed


\end{document}